\documentclass[12pt]{article}
\usepackage{lscape}
\usepackage{pdflscape}
\usepackage{amssymb}
\usepackage{amsfonts}
\usepackage{amsmath}
\usepackage{amsbsy}
\newcommand{\be}{\begin{equation}}
\newcommand{\ee}{\end{equation}}
\newcommand{\bea}{\begin{eqnarray}}
\newcommand{\eea}{\end{eqnarray}}
\newcommand{\ba}{\begin{array}}
\newcommand{\ea}{\end{array}}

\newcommand{\bc}{\begin{center}}
\newcommand{\ec}{\end{center}}
\newcommand{\ben}{\begin{enumerate}}
\newcommand{\een}{\end{enumerate}}
\newcommand{\bfi}{\begin{figure}}
\newcommand{\efi}{\end{figure}}

\newcommand{\bq}{\begin{quote}}
\newcommand{\eq}{\end{quote}}
\newcommand{\bqu}{\begin{quotation}}
\newcommand{\equ}{\end{quotation}}
\newenvironment{emphit}{\begin{itemize}}{\end{itemize}}
\newcommand{\bemp}{\begin{emphit}}
\newcommand{\eemp}{\end{emphit}}

\newcommand{\bt}{\begin{tabular}}
\newcommand{\et}{\end{tabular}}

\newtheorem{myth}{Theorem}[section]
\newtheorem{mylem}{Lemma}[section]

\newtheorem{mydef}{Definition}[section]
\newtheorem{myrem}{Remark}[section]
\newtheorem{myexam}{Example}[section]

\begin{document}
\date{}
\title{Varieties of groupoids and quasigroups generated by linear-bivariate polynomials over the ring $\mathbb{Z}_n$\footnote{2010 Mathematics Subject
Classification. Primary 20N02, 20NO5; Secondary 05B15}
\thanks{{\bf Keywords and Phrases :} groupoids, quasigroups, linear-bivariate polynomials}}
\author{E. Ilojide \\
Department of Mathematics,\\
Obafemi Awolowo University,\\
Ile-Ife 220005, Nigeria.\\
emmailojide@yahoo.com\and
T. G. Jaiy\'e\d ol\'a\thanks{All correspondence to be addressed to this author.} \\
Department of Mathematics,\\
Obafemi Awolowo University,\\
Ile-Ife 220005, Nigeria.\\
jaiyeolatemitope@yahoo.com\\tjayeola@oauife.edu.ng\and
O. O. Owojori\\
Department of Mathematics,\\
Obafemi Awolowo University,\\
Ile-Ife 220005, Nigeria.\\
walejori@oauife.edu.ng,walejori@yahoo.com}\maketitle
\begin{abstract}
Some varieties of groupoids and quasigroups generated by
linear-bivariate polynomials $P(x,y)=a+bx+cy$ over the ring
$\mathbb{Z}_n$ are studied. Necessary and sufficient conditions for
such groupoids and quasigroups to obey identities which involve one,
two, three (e.g. Bol-Moufang type) and four variables w.r.t. $a$, $b$ and $c$ are established.
Necessary and sufficient conditions for such groupoids and
quasigroups to obey some inverse properties w.r.t. $a$, $b$ and $c$ are also established.
This class of
groupoids and quasigroups are found to belong to some varieties of
groupoids and quasigroups such as medial groupoid(quasigroup),
F-quasigroup, semi automorphic inverse property groupoid(quasigroup) and automorphic inverse property
groupoid(quasigroup).
\end{abstract}
\section{Introduction}

\subsection{Groupoids, Quasigroups and Identities}
Let $G$  be a non-empty set. Define a binary operation ($\cdot$) on
$G$. $(G,\cdot)$ is called a groupoid if $G$ is closed under the
binary operation ($\cdot$). A groupoid $(G,\cdot)$ is called a
quasigroup if the equations $a\cdot x=b$ and $y\cdot c=d$ have
unique solutions for $x$ and $y$ for all $a,b,c,d\in G$. A
quasigroup $(G,\cdot)$ is called a loop if there exists a unique
element $e\in G$ called the identity element such that $x\cdot
e=e\cdot x=x$ for all $x\in G$.

A function $f~:~S\times S\to S$ on a finite set $S$ of size $n>0$ is
said to be a Latin square (of order $n$) if for any value $a\in S$
both functions $f(a,\cdot )$ and $f(\cdot , a)$ are permutations of
$S$. That is, a Latin square is a square matrix with $n^{2}$ entries
of $n$ different elements, none of them occurring more than once
within any row or column of the matrix.

\begin{mydef}
A pair of Latin squares $f_1(\cdot ,\cdot )$ and $f_2(\cdot ,\cdot
)$ is said to be orthogonal if the pairs $\big(f_1(x ,y),f_2(x,y)\big)$ are all distinct, as $x$ and $y$ vary.
\end{mydef}
\paragraph{}
For associative binary systems, the concept of an inverse element is
only meaningful if the system has an identity element. For example,
in a group $(G,\cdot )$ with identity element $e\in G$, if $x\in G$
then the inverse element for $x$ is the element $x^{-1}\in G$ such
that
\begin{displaymath}
x\cdot x^{-1}=x^{-1}\cdot x=e.
\end{displaymath}
In a loop $(G,\cdot )$ with identity element $e$, the left inverse
element of $x\in G$ is the element $x^\lambda\in G$ such that
\begin{displaymath}
x^\lambda\cdot x=e
\end{displaymath}
while the right inverse element of $x\in G$ is the element
$x^\rho\in G$ such that
\begin{displaymath}
x\cdot x^\rho=e
\end{displaymath}
In case $(G,\cdot )$ is a quasigroup, then for each $x\in G$, the
elements $x^\rho\in G$ and $x^\lambda\in G$ such that
$xx^\rho=e^\rho$ and $x^\lambda x=e^\lambda$ are called the right
and left inverse elements of $x$ respectively. Here, $e^\rho\in G$
and $e^\lambda\in G$ satisfy the relations $xe^\rho =x$ and
$e^\lambda x=x$ for all $x\in G$ and are respectively called the
right and left identity elements. Whenever $e^\rho =e^\lambda$, then
$(G,\cdot )$ becomes a loop.

In case $(G,\cdot )$ is a groupoid, then for each $x\in G$, the
elements $x^\rho\in G$ and $x^\lambda\in G$ such that
$xx^\rho=e_\rho(x)$ and $x^\lambda x=e_\lambda (x)$ are called the
right and left inverse elements of $x$ respectively. Here,
$e_\rho(x)\in G$ and $e_\lambda (x)\in G$ satisfy the relations
$xe_\rho(x) =x$ and $e_\lambda (x)x=x$ for each $x\in G$ and are
respectively called the local right and local left identity elements
of $x$. Whenever $e_\rho (x) =e_\lambda (x)$, then we simply write
$e(x)=e_\rho (x) =e_\lambda (x)$ and call it the local identity of
$x$.

The basic text books on quasigroups, loops are Pflugfelder
\cite{immaref:1}, Bruck \cite{immaref:48}, Chein, Pflugfelder and
Smith \cite{immaref:11}, Dene and Keedwell \cite{immaref:2},
Goodaire, Jespers and Milies \cite{phd42}, Sabinin
\cite{immaref:13}, Smith \cite{immaref:14}, Ja\'iy\'e\d ol\'a
\cite{immaref:12} and Vasantha Kandasamy \cite{phd75}.
\paragraph{}
Groupoids, quasigroups and loops are usually studied relative to
properties or identities. If a groupoid, quasigroup or loop obeys a
particular identity, then such types of groupoids, quasigroups or
loops are said to form a variety. In this work, our focus will be on
groupoids and quasigroups. Some identities that describe groupoids
and quasigroups which would be of interest to us here are
categorized as follows:
\begin{description}
\item[(A)] Those identities which involve one element only on each side of the equality sign:
\begin{equation}\label{imeq:1}
aa=a\qquad\qquad\textrm{idempotent law}
\end{equation}
\begin{equation}\label{imeq:2}
aa = bb\qquad\qquad\textrm{unipotent law}
\end{equation}
\item[(B)] Those identities which involve two elements on one or both sides of the equality sign:
\begin{equation}\label{imeq:3}
ab = ba\qquad\qquad\textrm{commutative law}
\end{equation}
\begin{equation}\label{imeq:4}
(ab)b = a\qquad\qquad\textrm{Sade right Keys law}
\end{equation}
\begin{equation}\label{imeq:5}
b(ba) = a\qquad\qquad\textrm{Sade left keys law}
\end{equation}
\begin{equation}\label{imeq:6}
(ab)b = a(bb)\qquad\qquad\textrm{right alternative  law}
\end{equation}
\begin{equation}\label{imeq:7}
b(ba) = (bb)a\qquad\qquad\textrm{left alternative  law}
\end{equation}
\begin{equation}\label{imeq:8}
a(ba) = (ab)a\qquad\qquad\textrm{medial alternative  law}
\end{equation}
\begin{equation}\label{imeq:9}
a(ba) = b\qquad\qquad\textrm{law of right semisymmetry}
\end{equation}
\begin{equation}\label{imeq:10}
(ab)a = b\qquad\qquad\textrm{law of left semisymmetry}
\end{equation}
\begin{equation}\label{imeq:11}
a(ab) = ba\qquad\qquad\textrm{Stein first law }
\end{equation}
\begin{equation}\label{imeq:12}
a(ba) = (ba)a\qquad\qquad\textrm{Stein second law }
\end{equation}
\begin{equation}\label{imeq:13}
a(ab) = (ab)b\qquad\qquad\textrm{Schroder first  law }
\end{equation}
\begin{equation}\label{imeq:14}
(ab) (ba) = a\qquad\qquad\textrm{Schroder second  law }
\end{equation}
\begin{equation}\label{imeq:15}
(ab) (ba) = b\qquad\qquad\textrm{Stein third law }
\end{equation}
\begin{equation}\label{imeq:16}
ab = a\qquad\qquad\textrm{Sade right translation law }
\end{equation}
\begin{equation}\label{imeq:17}
ab = b\qquad\qquad\textrm{Sade left translation law }
\end{equation}
\item[(C)] Those identities which involve three distinct elements on one or both sides of the equality sign:
\begin{equation}\label{imeq:18}
(ab)c = a(bc)\qquad\qquad\textrm{associative law }
\end{equation}
\begin{equation}\label{imeq:19}
a(bc) = c(ab)\qquad\qquad\textrm{law of cyclic associativity }
\end{equation}
\begin{equation}\label{imeq:20}
(ab)c = (ac)b\qquad\qquad\textrm{law of right permutability }
\end{equation}
\begin{equation}\label{imeq:21}
a(bc) = b(ac)\qquad\qquad\textrm{law of left permutability }
\end{equation}
\begin{equation}\label{imeq:22}
a(bc) = c(ba)\qquad\qquad\textrm{Abel-Grassman law }
\end{equation}
\begin{equation}\label{imeq:23}
(ab)c = a(cb)\qquad\qquad\textrm{commuting product law }
\end{equation}
\begin{equation}\label{imeq:24}
c(ba) = (bc)a\qquad\qquad\textrm{dual of commuting product }
\end{equation}
\begin{equation}\label{imeq:25}
(ab) (bc) = ac\qquad\qquad\textrm{Stein fourth law }
\end{equation}
\begin{equation}\label{imeq:26}
(ba) (ca) = bc\qquad\qquad\textrm{law of right transitivity }
\end{equation}
\begin{equation}\label{imeq:27}
(ab) (ac) = bc\qquad\qquad\textrm{law of left transitivity }
\end{equation}
\begin{equation}\label{imeq:28}
(ab) (ac) = cb\qquad\qquad\textrm{Schweitzer law }
\end{equation}
\begin{equation}\label{imeq:29}
(ba) (ca) = cb\qquad\qquad\textrm{dual of Schweitzer law }
\end{equation}
\begin{equation}\label{imeq:30}
(ab)c = (ac) (bc)\qquad\qquad\textrm{law of right self-distributivity law }
\end{equation}
\begin{equation}\label{imeq:31}
c(ba) = (cb) (ca)\qquad\qquad\textrm{law of left self-distributivity law }
\end{equation}
\begin{equation}\label{imeq:32}
(ab)c = (ca) (bc)\qquad\qquad\textrm{law of right abelian distributivity }
\end{equation}
\begin{equation}\label{imeq:33}
c(ba) = (cb) (ac)\qquad\qquad\textrm{law of left abelian distributivity }
\end{equation}
\begin{equation}\label{imeq:34}
(ab) (ca) = [a(bc)]a\qquad\qquad\textrm{Bruck-Moufang identity }
\end{equation}
\begin{equation}\label{imeq:35}
(ab) (ca) = a[bc)a]\qquad\qquad\textrm{dual of Bruck-Moufang identity }
\end{equation}
\begin{equation}\label{imeq:36}
[(ab)c]b = a[b(cb)]\qquad\qquad\textrm{Moufang identity }
\end{equation}
\begin{equation}\label{imeq:37}
[(bc)b]a = b[c(ba)]\qquad\qquad\textrm{Moufang identity }
\end{equation}
\begin{equation}\label{imeq:38}
[(ab)c]b = a[(bc)b]\qquad\qquad\textrm{right Bol identity }
\end{equation}
\begin{equation}\label{imeq:39}
[b(cb)]a = b[c(ba)]\qquad\qquad\textrm{left Bol identity }
\end{equation}
\begin{equation}\label{imeq:40}
[(ab)c]a = a[b(ca)]\qquad\qquad\textrm{extra law }
\end{equation}
\begin{equation}\label{imeq:40.1}
[(ba)a]c = b[(aa)c]\qquad\qquad\textrm{$\textrm{RC}_4$ law}
\end{equation}
\begin{equation}\label{imeq:40.2}
[b(aa)]c = b[a(ac)]\qquad\qquad\textrm{$\textrm{LC}_4$ law }
\end{equation}
\begin{equation}\label{imeq:40.3}
(aa)(bc) = [a(ab)]c\qquad\qquad\textrm{$\textrm{LC}_2$ law }
\end{equation}
\begin{equation}\label{imeq:40.4}
[(bc)a]a = b[(ca)a]\qquad\qquad\textrm{$\textrm{RC}_1$ law }
\end{equation}
\begin{equation}\label{imeq:40.5}
[a(ab)]c = a[a(bc)]\qquad\qquad\textrm{$\textrm{LC}_1$ law }
\end{equation}
\begin{equation}\label{imeq:40.6}
(bc)(aa) = b[(ca)a]\qquad\qquad\textrm{$\textrm{RC}_2$ law }
\end{equation}
\begin{equation}\label{imeq:40.7}
[(aa)b]c = a[a(bc)]\qquad\qquad\textrm{$\textrm{LC}_3$ law }
\end{equation}
\begin{equation}\label{imeq:40.8}
[(bc)a]a = b[c(aa)]\qquad\qquad\textrm{$\textrm{RC}_3$ law }
\end{equation}
\begin{equation}\label{imeq:40.9}
[(ba)a]c = b[a(ac)]\qquad\qquad\textrm{C-law}
\end{equation}
\begin{equation}\label{imeq:41}
a[b(ca)] = cb\qquad\qquad\textrm{Tarski law }
\end{equation}
\begin{equation}\label{imeq:42}
a[(bc) (ba)] = c\qquad\qquad\textrm{ Neumann law }
\end{equation}
\begin{equation}\label{imeq:43}
(ab) (ca) = (ac) (ba)\qquad\qquad\textrm{ specialized medial law }
\end{equation}
\item[(D)] Those involving four elements:
\begin{equation}\label{imeq:44}
(ab) (cd) = (ad) (cb)\qquad\qquad\textrm{ first rectangle rule }
\end{equation}
\begin{equation}\label{imeq:45}
(ab) (ac) = (db) (dc)\qquad\qquad\textrm{ second rectangle rule }
\end{equation}
\begin{equation}\label{imeq:46}
(ab) (cd) = (ac) (bd)\qquad\qquad\textrm{ internal mediality or
medial law }
\end{equation}
\item[(E)] Those involving left or right inverse elements:
\begin{equation}\label{imeq:44.1}
x^\lambda\cdot xy=y\qquad\qquad\textrm{left inverse property}
\end{equation}
\begin{equation}\label{imeq:45.1}
yx\cdot x^\rho=y\qquad\qquad\textrm{right inverse property}
\end{equation}
\begin{equation}\label{imeq:47}
x(yx)^\rho=y^\rho~\textrm{or}~(xy)^\lambda
x=y^\lambda\qquad\qquad\textrm{weak inverse property(WIP)}
\end{equation}
\begin{equation}\label{imeq:48}
xy\cdot x^\rho =y~\textrm{or}~ x\cdot yx^\rho
=y~\textrm{or}~x^\lambda\cdot (yx)=y~\textrm{or}~x^\lambda y\cdot
x=y~\textrm{cross inverse property(CIP)}
\end{equation}
\begin{equation}\label{imeq:49}
(xy)^\rho=x^\rho y^\rho~\textrm{or}~(xy)^\lambda =x^\lambda
y^\lambda~\textrm{automorphic inverse property (AIP)}
\end{equation}
\begin{equation}\label{imeq:50}
(xy)^\rho= y^\rho x^\rho~\textrm{or}~(xy)^\lambda =y^\lambda
x^\lambda~\textrm{anti-automorphic inverse property (AAIP) }
\end{equation}
\begin{equation}\label{imeq:51}
(xy\cdot x)^\rho= x^\rho y^\rho\cdot x^\rho~\textrm{or}~(xy\cdot
x)^\lambda =x^\lambda y^\lambda\cdot x^\lambda
~\textrm{semi-automorphic inverse property (SAIP) }
\end{equation}
\end{description}

\begin{mydef}(Trimedial Quasigroup)

A quasigroup is trimedial if every subquasigroup generated by three
elements is medial.
\end{mydef}

Medial quasigroups have also been called abelian, entropic, and other names,
while trimedial quasigroups have also been called triabelian, terentropic, etc.

There are two distinct, but related, generalizations of trimedial quasigroups.
The variety of semimedial quasigroups(also known as weakly abelian, weakly medial, etc.)
is defined by the equations
\begin{equation}\label{imeq:54}
xx\cdot yz = xy\cdot xz
\end{equation}
\begin{equation}\label{imeq:55}
zy\cdot xx = zx\cdot yz
\end{equation}

\begin{mydef}(Semimedial Quasigroup)

A quasigroup satisfying (\ref{imeq:54}) (resp. (\ref{imeq:55})) is
said to be left (resp. right) semimedial.

\end{mydef}

\begin{mydef}(Medial-Like Identities)

A groupoid or quasigroup is called an external medial groupoid or
quasigroup if it obeys the identity
\begin{equation}\label{imeq:56}
ab\cdot cd = db\cdot ca\qquad\qquad\textrm{external medial or
paramediality law}
\end{equation}
A groupoid or quasigroup is called a palindromic groupoid or
quasigroup if it obeys the identity
\begin{equation}\label{imeq:56.1}
ab\cdot cd = dc\cdot ba\qquad\qquad\textrm{palidromity law}
\end{equation}
Other medial like identities of the form $(ab) (cd)=(\pi (a)\pi (b))
(\pi (c)\pi (d))$, where $\pi$ is a certain permutation on
$\{a,b,c,d\}$ are given as follows:

\begin{equation}\label{imeq:56.2}
ab\cdot cd=ab\cdot dc\qquad\qquad\textrm{C$_1$}\end{equation}
\begin{equation}\label{imeq:56.3}
ab\cdot cd=ba\cdot cd\qquad\qquad\textrm{C$_2$}\end{equation}
\begin{equation}\label{imeq:56.4}
ab\cdot cd=ba\cdot dc\qquad\qquad\textrm{C$_3$}\end{equation}
\begin{equation}\label{imeq:56.5}
ab\cdot cd=cd\cdot ab\qquad\qquad\textrm{C$_4$}\end{equation}
\begin{equation}\label{imeq:56.6}
ab\cdot cd=cd\cdot ba\qquad\qquad\textrm{C$_5$}\end{equation}
\begin{equation}\label{imeq:56.7}
ab\cdot cd=dc\cdot ab\qquad\qquad\textrm{C$_6$}\end{equation}
\begin{equation}\label{imeq:56.8}
ab\cdot cd=ac\cdot db\qquad\qquad\textrm{CM$_1$}\end{equation}
\begin{equation}\label{imeq:56.9}
ab\cdot cd=ad\cdot bc\qquad\qquad\textrm{CM$_2$}\end{equation}
\begin{equation}\label{imeq:56.10}
ab\cdot cd=ad\cdot cb\qquad\qquad\textrm{CM$_3$}\end{equation}
\begin{equation}\label{imeq:56.11}
ab\cdot cd=bc\cdot ad\qquad\qquad\textrm{CM$_4$}\end{equation}
\begin{equation}\label{imeq:56.12}
ab\cdot cd=bc\cdot da\qquad\qquad\textrm{CM$_5$}\end{equation}
\begin{equation}\label{imeq:56.13}
ab\cdot cd=bd\cdot ac\qquad\qquad\textrm{CM$_6$}\end{equation}
\begin{equation}\label{imeq:56.14}
ab\cdot cd=bd\cdot ca\qquad\qquad\textrm{CM$_7$}\end{equation}
\begin{equation}\label{imeq:56.15}
ab\cdot cd=ca\cdot bd\qquad\qquad\textrm{CM$_8$}\end{equation}
\begin{equation}\label{imeq:56.16}
ab\cdot cd=ca\cdot db\qquad\qquad\textrm{CM$_9$}\end{equation}
\begin{equation}\label{imeq:56.17}
ab\cdot cd=cb\cdot ad\qquad\qquad\textrm{CM$_{10}$}\end{equation}
\begin{equation}\label{imeq:56.18}
ab\cdot cd=cb\cdot da\qquad\qquad\textrm{CM$_{11}$}\end{equation}
\begin{equation}\label{imeq:56.19}
ab\cdot cd=da\cdot bc\qquad\qquad\textrm{CM$_{12}$}\end{equation}
\begin{equation}\label{imeq:56.20}
ab\cdot cd=da\cdot cb\qquad\qquad\textrm{CM$_{13}$}\end{equation}
\begin{equation}\label{imeq:56.21}
ab\cdot cd=db\cdot ac\qquad\qquad\textrm{CM$_{14}$}\end{equation}
\end{mydef}

The variety of F-quasigroups was introduced by Murdoch
\cite{immaref:27}.
\begin{mydef}(F-quasigroup)

An F-quasigroup is a quasigroup that obeys the identities
\begin{equation}\label{imeq:57}
x\cdot yz = xy\cdot (x\backslash x)z\qquad\qquad\textrm{left F-law}
\end{equation}
\begin{equation}\label{imeq:58}
zy\cdot x = z(x/x)\cdot yx\qquad\qquad\textrm{right F-law}
\end{equation}
A quasigroup satisfying (\ref{imeq:57}) (resp. (\ref{imeq:58})) is
called a left (resp. right) F-quasigroup.
\end{mydef}

\begin{mydef}(E-quasigroup)

An E-quasigroup is a quasigroup that obeys the identities
\begin{equation}\label{imeq:57.1}
x\cdot yz=e_\lambda (x)y\cdot xz\qquad\qquad\textrm{E}_\textrm{l}~\textrm{law}
\end{equation}
\begin{equation}\label{imeq:58.1}
zy\cdot x=zx\cdot y e_\rho (x)\qquad\qquad\textrm{E}_\textrm{r}~\textrm{Law}
\end{equation}
A quasigroup satisfying (\ref{imeq:57.1}) (resp. (\ref{imeq:58.1})) is
called a left (resp. right) E-quasigroup.
\end{mydef}
Some identities will make a quasigroup to be a loop, such are
discussed in Keedwell \cite{immaref:28,immaref:29}.

\begin{mydef}(Linear Quasigroup and T-quasigroup)
A quasigroup $(Q,\cdot )$ of the form $x\cdot y = x\alpha+y\beta+c$ where $(Q, +)$ is a group, $\alpha$ is its automorphism and
$\beta$ is a permutation of the set $Q$, is called a left linear quasigroup.

A quasigroup $(Q,\cdot )$ of the form $x\cdot y = x\alpha+y\beta+c$ where $(Q, +)$ is a group, $\beta$ is its automorphism and
$\alpha$ is a permutation of the set $Q$, is called a right linear quasigroup.

A T-quasigroup is a quasigroup $(Q, \cdot)$ defined over an abelian group $(Q, +)$ by $x \cdot y=c+x\alpha +y\beta$,
where $c$ is a fixed element of $Q$ and $\alpha$ and $\beta$ are both automorphisms of the group $(Q, +)$.
\end{mydef}
\paragraph{}Whenever one considers mathematical objects defined in some abstract
manner, it is usually desirable to determine that such objects
exist. Although occasionally this is accomplished by means of an
abstract existential argument, most frequently, it is carried out
through the presentation of a suitable example, often one which has
been specifically constructed for the purpose. An example is the
solution to the open problem of the axiomization of rectangular
quasigroups and loops by Kinyon and Phillips \cite{phd125} and the
axiomization of trimedial quasigroups by Kinyon and Phillips
\cite{immaref:25,immaref:26}.

Chein et. al. \cite{immaref:11} presents a survey of various methods
of construction which has been used in the literature to generate
examples of groupoids and quasigroups. Many of these constructions
are ad hoc-designed specifically to produce a particular example;
while others are of more general applicability. More can be found on
the construction of $(r,s,t)$-inverse quasigroups in Keedwell and
Shcherbacov \cite{immaref:21,immaref:22}, idempotent medial
quasigroups in Kr\v{c}adinac and Volenec \cite{immaref:31} and
quasigroups of Bol-Moufang type in Kunen
\cite{immaref:23,immaref:24}.
\begin{myrem}
In the survey of methods of construction of varieties and types of
quasigroups highlighted in Chein et. al. \cite{immaref:11}, it will
be observed that some other important types of quasigroups that obey
identities (\ref{imeq:1}) to (\ref{imeq:58.1}) are not mentioned.
Also, examples of methods of construction of such varieties that are
groupoids are also scarce or probably not in existence by our
search. In Theorem 1.4 of Kirnasovsky \cite{kir}, the author characterized T-quasigroups with a score and two identities from among identities (\ref{imeq:1}) to (\ref{imeq:58.1}). The present work thus proves some results with which such
groupoids and quasigroups can be constructed.
\end{myrem}

\subsection{Univariate and Bivariate Polynomials}
Consider the following definitions.
\begin{mydef}\label{imdef:10}

A polynomial $P(x) = a_{0} +a_{1}x + \cdots +a_{n}x^{n}$,
$n\in\mathbb{N}$ is said to be a permutation polynomial over a
finite ring $R$ if the mapping defined by $P$ is a bijection on $R$.
\end{mydef}
\begin{mydef}\label{imdef:11}

A bivariate polynomial is a polynomial in two variables, $x$  and
$y$ of the form $P(x,y) =
\displaystyle\Sigma_{i,j}a_{ij}x^{i}y^{j}$.
\end{mydef}
\begin{mydef}(Bivariate Polynomial Representing a Latin Square)\label{imdef:11.1}

A bivariate polynomial $P(x,y)$ over $\mathbb{Z}_n$ is said to
represent (or generate) a Latin square if $(\mathbb{Z}_n,\ast )$ is
a quasigroup where $\ast~:~\mathbb{Z}_n\times \mathbb{Z}_n\to
\mathbb{Z}_n$ is defined by $x\ast y=P(x,y)$ for all $x,y\in
\mathbb{Z}_n$.
\end{mydef}

Mollin and Small \cite{immaref:16} considered the problem
of characterizing permutation polynomials. They established
conditions on the coefficients of a polynomial which are necessary
and sufficient for it to represent a permutation.

Shortly after, Rudolf and Mullen \cite{immaref:17} provided a brief
survey of the main known classes of permutation polynomials over a
finite field and discussed some problems concerning permutation
polynomials (PPs). They described several applications of
permutations which indicated why the study of permutations is of
interest. Permutations of finite fields have become of considerable
interest in the construction of cryptographic systems for the secure
transmission of data. Thereafter, the same authors in their
paper \cite{immaref:18}, described some results that had appeared
after their earlier work including two major breakthroughs.

Rivest \cite{immaref:19} studied permutation polynomials
over the ring $(\mathbb{Z}_n,+,\cdot )$ where $n$ is a power of $2$:
$n=2^w$. This is based on the fact that modern computers perform
computations modulo $2^w$ efficiently (where $w=2,8,16,32$ or $64$
is the word size of the machine), and so it was of interest to study
PPs modulo a power of $2$. Below is an important result from his
work which is relevant to the present study.

\begin{myth}\label{immathm:7}(Rivest \cite{immaref:19})

A bivariate polynomial $P(x,y) =
\displaystyle\Sigma_{i,j}a_{ij}x^{i}y^{j}$ represents a Latin square
modulo $n = 2^{w}$, where $w\geq 2$, if and only if the four
univariate polynomials $P(x,0)$, $P(x,1)$, $P(0,y)$, and  $P(1,y)$
are all permutation polynomial modulo $n$.
\end{myth}

Vadiraja and Shankar \cite{immaref:3} motivated by the work
of Rivest continued the study of permutation polynomials over the
ring $(\mathbb{Z}_n,+,\cdot )$ by studying Latin squares represented
by linear and quadratic bivariate polynomials over $\mathbb{Z}_n$
when $n\ne 2^w$ with the characterization of some PPs. Some of the
main results they got are stated below.

\begin{myth}\label{immathm:9}(Vadiraja and Shankar \cite{immaref:3})

A bivariate linear polynomial $a + bx + cy$ represents a Latin
square over $\mathbb{Z}_{n}$, $n\ne 2^w$ if and only if one of the
following equivalent conditions is satisfied:
\begin{description}
\item[(i)] both b and c are coprime with $n$;
\item[(ii)] $a + bx$, $a + cy$, $(a + c) + bx$ and $(a + b) + cy$ are all
permutation polynomials modulo $n$.
\end{description}
\end{myth}
\begin{myrem}
It must be noted that $P(x,y)=a + bx + cy$ represents a groupoid over $\mathbb{Z}_{n}$. $P(x,y)$ represents a quasigroup over $\mathbb{Z}_{n}$ if and only if $(\mathbb{Z}_{n},P)$ is a T-quasigroup. Hence whenever $(\mathbb{Z}_{n},P)$ is a groupoid and not a quasigroup, $(\mathbb{Z}_{n},P)$ is neither a T-quasigroup nor left linear quasigroup nor right linear quasigroup. Thus, the present study considers both T-quasigroup and non-T-quasigroup.
\end{myrem}
\begin{myth}\label{immathm:10}(Vadiraja and Shankar \cite{immaref:3})

If $P(x,y)$ is a bivariate polynomial having no cross term, then
$P(x,y)$ gives a Latin square if and only if $P(x,0)$ and $P(0,y)$
are permutation polynomials.
\end{myth}

The authors were able to establish the fact that Rivest's result for
a bivariate polynomial over $\mathbb{Z}_n$ when $n=2^w$ is true for
a linear-bivariate polynomial over $\mathbb{Z}_n$ when $n\ne 2^w$.
Although the result of Rivest was found not to be true for
quadratic-bivariate polynomials over $\mathbb{Z}_n$ when $n\ne 2^w$
with the help of counter examples, nevertheless some of such squares
can be forced to be Latin squares by deleting some equal numbers of
rows and columns.

Furthermore, Vadiraja and Shankar \cite{immaref:3} were able to find
examples of pairs of orthogonal Latin squares generated by bivariate
polynomials over $\mathbb{Z}_n$ when $n\ne 2^w$ which was found
impossible by Rivest for bivariate polynomials over $\mathbb{Z}_n$
when $n= 2^w$.

\subsection{Some Important Results on Medial-Like Identities}
Some important results which we would find useful in our study are
stated below.
\begin{myth}\label{immathm:13}(Polonijo \cite{immaref:30})

For any groupoid $(Q,\cdot)$, any two of the three identities
(\ref{imeq:46}), (\ref{imeq:56}) and (\ref{imeq:56.1}) imply the
third one.
\end{myth}
\begin{myth}\label{immathm:14}(Polonijo \cite{immaref:30})

Let $(Q,\cdot)$ be a commutative groupoid. Then $(Q,\cdot)$ is
palindromic. Furthermore, the constraints (\ref{imeq:46}) and
(\ref{imeq:56}) are equivalent, i.e a commutative groupoid
$(Q,\cdot)$ is internally medial if and only if it is externally
medial.
\end{myth}
\begin{myth}\label{immathm:24}(Polonijo \cite{immaref:30})

For any quasigroup $(Q,\cdot)$ and $i\in \{1, 2,..., 6\}$, C$_{i}$
is valid if and only if the quasigroup is commutative.
\end{myth}
\begin{myth}\label{immathm:25}(Polonijo \cite{immaref:30})

For any quasigroup $(Q,\cdot)$ and $i\in \{1, 2,..., 14\}$, CM$_{i}$
holds if and only if the quasigroup is both commutative and
internally medial.
\end{myth}
\begin{myth}\label{immathm:26}(Polonijo \cite{immaref:30})

For any quasigroup $(Q,\cdot)$ and $i\in \{1, 2,..., 14\}$, CM$_{i}$
is valid if and only if the quasigroup is both commutative and
externally medial.
\end{myth}

\begin{myth}\label{prop}
A quasigroup $(Q,\cdot )$ is palindromic if and only if there exists an automorphism $\alpha$ such that
\begin{displaymath}
\alpha (x\cdot y)=y\cdot x~\forall~x,y\in Q
\end{displaymath}
holds.
\end{myth}
\paragraph{}It is important to study the characterization of varieties of
groupoids and quasigroups represented by linear-bivariate
polynomials over the ring $\mathbb{Z}_n$ even though very few of
such have been sighted as examples in the past.

\section{Main Results}
\begin{myth}\label{immathm:30}

Let $P(x,y)=a+bx+cy$ be a linear bivariate polynomial over
$\{\mathbb{Z}_n,\mathbb{Z}_p\}$ such that "HYPO" is true. $P(x,y)$ represents a "NAME" $\{$groupoid, quasigroup$\}$
$\{(\mathbb{Z}_n,P),(\mathbb{Z}_p,P)\}$ over $\{\mathbb{Z}_n,\mathbb{Z}_p\}$ if and only if "N and S" is true. (Table~\ref{osbornloop9})
\end{myth}
{\bf Proof}\\
There are 66 identities for which the theorem above is true for in a groupoid or quasigroup. For the sake of space, we shall only demonstrate the proof for one identity for each category.
\begin{description}
\item[(A)] Those identities which involve one element only on each side of the equality sign:
\begin{mylem}\label{immalem:7}
Let $P(x,y)=a+bx+cy$ be a linear bivariate polynomial over
$\mathbb{Z}_n$. $P(x,y)$ represents a unipotent groupoid
$(\mathbb{Z}_n,P)$ over $\mathbb{Z}_n$ if and only if $(b+c)(x-y)=0$
for all $x,y\in \mathbb{Z}_n$.
\end{mylem}
{\bf Proof}\\
$P(x,y)$ satisfies the unipotent law $\Leftrightarrow$ $P(x,x)=P(y,y)$
$\Leftrightarrow$ $a+bx+cx=a+by+cy$ $\Leftrightarrow$ $a+bx-cx-a-by-cy=0$
$\Leftrightarrow$ $(b+c)(x-y)=0$ as required.
\begin{mylem}\label{immalem:7.1}
Let $P(x,y)=a+bx+cy$ be a linear bivariate polynomial over
$\mathbb{Z}_n$. $P(x,y)$ represents a unipotent quasigroup
$(\mathbb{Z}_n,P)$ over $\mathbb{Z}_n$ if and only if $(b+c)(x-y)=0$ and $(b,n)=(c,n)=1$
for all $x,y\in \mathbb{Z}_n$.
\end{mylem}
{\bf Proof}\\
This is proved by using Lemma~\ref{immalem:7} and Theorem~\ref{immathm:9}
\begin{myth}\label{immathm:31}
Let $P(x,y)=a+bx+cy$ be a linear bivariate polynomial over
$\mathbb{Z}_n$. $P(x,y)$ represents a unipotent groupoid
$(\mathbb{Z}_n,P)$ over $\mathbb{Z}_n$ if and only if $b+c\equiv 0(\bmod{n})$.
\end{myth}
{\bf Proof}\\
This is proved by using Lemma~\ref{immalem:7}.
\begin{myth}\label{immathm:31.1}
Let $P(x,y)=a+bx+cy$ be a linear bivariate polynomial over
$\mathbb{Z}_n$. $P(x,y)$ represents a unipotent quasigroupp
$(\mathbb{Z}_n,P)$ over $\mathbb{Z}_n$ if and only if $b+c\equiv 0(\bmod{n})$ and $(b,n)=(c,n)=1$.
\end{myth}
{\bf Proof}\\
This is proved by using Lemma~\ref{immalem:7.1}.
\begin{myexam}
$P(x,y)=5x+y$ is a linear bivariate polynomial over $\mathbb{Z}_6$.
$(\mathbb{Z}_6,P)$ is a unipotent groupoid over $\mathbb{Z}_6$.
\end{myexam}
\begin{myexam}
$P(x,y)=1+5x+y$ is a linear bivariate polynomial over $\mathbb{Z}_6$.
$(\mathbb{Z}_6,P)$ is a unipotent quasigroup over $\mathbb{Z}_6$.
\end{myexam}
\item[(B)] Those identities which involve two elements on one or both sides of the equality sign:
\begin{mylem}\label{immalem:12.14}
Let $P(x,y)=a+bx+cy$ be a linear bivariate polynomial over
$\mathbb{Z}_n$. $P(x,y)$ represents a Stein third groupoid
$(\mathbb{Z}_n,P)$ over $\mathbb{Z}_n$ if and only if $a(1+b+c)+x(b^2+c^2)+y(2bc-1)=0$
for all $x,y\in \mathbb{Z}_n$.
\end{mylem}
{\bf Proof}\\
$P(x,y)$ satisfies the Stein third law $\Leftrightarrow P[P(x,y),P(y,x)]=y$
$\Leftrightarrow$ $a(1+b+c)+x(b^2+c^2)+y(2bc-1)=0$ as required.
\begin{mylem}\label{immalem:12.15}
Let $P(x,y)=a+bx+cy$ be a linear bivariate polynomial over
$\mathbb{Z}_n$. $P(x,y)$ represents a Stein third quasigroup
$(\mathbb{Z}_n,P)$ over $\mathbb{Z}_n$ if and only if $a(1+b+c)+x(b^2+c^2)+y(2bc-1)=0$
and $(b,n)=(c,n)=1$ for all $x,y\in \mathbb{Z}_n$.
\end{mylem}
{\bf Proof}\\
This is proved by using Lemma~\ref{immalem:12.14} and Theorem~\ref{immathm:9}
\begin{myth}\label{immathm:36.22}
Let $P(x,y)=a+bx+cy$ be a linear bivariate polynomial over
$\mathbb{Z}_n$. $P(x,y)$ represents a Stein third groupoid
$(\mathbb{Z}_n,P)$ over $\mathbb{Z}_n$ if and only if $b^2+c^2\equiv 0(\bmod{n})$, $2bc\equiv 1(\bmod{n})$
and $a=0$.
\end{myth}
{\bf Proof}\\
This is proved by using Lemma~\ref{immalem:12.14}.
\begin{myth}\label{immathm:36.23}
Let $P(x,y)=a+bx+cy$ be a linear bivariate polynomial over
$\mathbb{Z}_n$. $P(x,y)$ represents a Stein third quasigroup
$(\mathbb{Z}_n,P)$ over $\mathbb{Z}_n$ if and only if $b^2+c^2\equiv 0(\bmod{n})$, $2bc\equiv 1(\bmod{n})$
and $a=0$.
\end{myth}
{\bf Proof}\\
This is proved by using Lemma~\ref{immalem:12.15}.
\begin{myth}\label{immathm:36.24}
Let $P(x,y)=a+bx+cy$ be a linear bivariate polynomial over
$\mathbb{Z}_p$ such that $a\not =0$. $P(x,y)$ represents a Stein third groupoid
$(\mathbb{Z}_p,P)$ over $\mathbb{Z}_p$ if and only if $b^2+c^2\equiv 0(\bmod{p})$
and $2bc\equiv 1(\bmod{p})$.
\end{myth}
{\bf Proof}\\
This is proved by using Lemma~\ref{immalem:12.14}.
\begin{myth}\label{immathm:36.25}
Let $P(x,y)=a+bx+cy$ be a linear bivariate polynomial over
$\mathbb{Z}_p$ such that $a\not =0$. $P(x,y)$ represents a Stein third quasigroup
$(\mathbb{Z}_p,P)$ over $\mathbb{Z}_p$ if and only if $b^2+c^2\equiv 0(\bmod{p})$
and $2bc\equiv 1(\bmod{p})$.
\end{myth}
{\bf Proof}\\
This is proved by using Lemma~\ref{immalem:12.15}.
\begin{myexam}
$P(x,y)=2x+3y$ is a linear bivariate polynomial over $\mathbb{Z}_5$.
$(\mathbb{Z}_5,P)$ is a Stein third groupoid over $\mathbb{Z}_5$.
\end{myexam}
\begin{myexam}
$P(x,y)=2x+3y$ is a linear bivariate polynomial over $\mathbb{Z}_5$.
$(\mathbb{Z}_5,P)$ is a Stein third quasigroup over $\mathbb{Z}_5$.
\end{myexam}
\item[(C)] Those identities which involve three distinct elements on one or both sides of the equality sign:
\begin{mylem}\label{immalem:17}
Let $P(x,y)=a+bx+cy$ be a linear bivariate polynomial over
$\mathbb{Z}_n$. $P(x,y)$ represents an Abel-Grassman groupoid
$(\mathbb{Z}_n,P)$ over $\mathbb{Z}_n$ if and only if $(x-z)(b-c^{2})=0$
for all $x,z\in \mathbb{Z}_n$.
\end{mylem}
{\bf Proof}\\
$P(x,y)$ satisfies the Abel-Grassman law $\Leftrightarrow$ $P[x,P(y,z)]=P[z,P(y,x)]$
$\Leftrightarrow$ $P(x,a+by+cz)=P(z,a+by+cx)$ $\Leftrightarrow$ $a+bx+c(a+by+cz)=a+bz+c(a+by+cx)$
$\Leftrightarrow$ $(x-z)(b-c^{2})=0$ as required.
\begin{mylem}\label{immalem:17.1}
Let $P(x,y)=a+bx+cy$ be a linear bivariate polynomial over
$\mathbb{Z}_n$. $P(x,y)$ represents an Abel-Grassman quasigroup
$(\mathbb{Z}_n,P)$ over $\mathbb{Z}_n$ if and only if $(x-z)(b-c^{2})=0$ and $(b,n)=(c,n)=1$.
for all $x,y,z\in \mathbb{Z}_n$.
\end{mylem}
{\bf Proof}\\
This is proved by using Lemma~\ref{immalem:17} and Theorem~\ref{immathm:9}
\begin{myth}\label{immathm:41}
Let $P(x,y)=a+bx+cy$ be a linear bivariate polynomial over
$\mathbb{Z}_n$. $P(x,y)$ represents an Abel-Grassman groupoid
$(\mathbb{Z}_n,P)$ over $\mathbb{Z}_n$ if and only if $c^{2}\equiv b(\bmod{n})$.
\end{myth}
{\bf Proof}\\
This is proved by using Lemma~\ref{immalem:17}.
\begin{myth}\label{immathm:41.1}
Let $P(x,y)=a+bx+cy$ be a linear bivariate polynomial over
$\mathbb{Z}_n$. $P(x,y)$ represents an Abel-Grassman quasigroup
$(\mathbb{Z}_n,P)$ over $\mathbb{Z}_n$ if and only if $c^{2}\equiv b(\bmod{n})$ and $(b,n)=(c,n)=1$.
\end{myth}
{\bf Proof}\\
This is proved by using Lemma~\ref{immalem:17.1}.
\begin{myexam}
$P(x,y)=2+4x+2y$ is a linear bivariate polynomial over $\mathbb{Z}_6$.
$(\mathbb{Z}_6,P)$ is an Abel-Grassman groupoid over $\mathbb{Z}_6$.
\end{myexam}
\begin{myexam}
$P(x,y)=2+4x+2y$ is a linear bivariate polynomial over $\mathbb{Z}_5$.
$(\mathbb{Z}_5,P)$ is an Abel-Grassman quasigroup over $\mathbb{Z}_5$.
\end{myexam}
\item[(D)] Those involving four elements:
\begin{mylem}\label{immalem:50}
Let $P(x,y)=a+bx+cy$ be a linear bivariate polynomial over
$\mathbb{Z}_n$. $P(x,y)$ represents an external medial groupoid
$(\mathbb{Z}_n,P)$ over $\mathbb{Z}_n$ if and only if $w(b^2-c^2)+z(c^2-b^2)=0$ for all $w,z\in \mathbb{Z}_n$.
\end{mylem}
{\bf Proof}\\
$P(x,y)$ satisfies the external medial law $\Leftrightarrow$ $P[P(w,x),P(y,z)]=P[P(z,x),P(y,w)]$
$\Leftrightarrow$ $a+b(a+bw+cx)+c(a+by+cz)=a+b(a+bz+cx)+c(a+by+cw)$ $\Leftrightarrow$ $w(b^2-c^2)+z(c^2-b^2)=0$ as required.
\begin{mylem}\label{immalem:50.1}
Let $P(x,y)=a+bx+cy$ be a linear bivariate polynomial over
$\mathbb{Z}_n$. $P(x,y)$ represents an external medial quasigroup
$(\mathbb{Z}_n,P)$ over $\mathbb{Z}_n$ if and only if $w(b^2-c^2)+z(c^2-b^2)=0$ and $(b,n)=(c,n)=1$
for all $w,z\in \mathbb{Z}_n$.
\end{mylem}
{\bf Proof}\\
This is proved by using Lemma~\ref{immalem:50} and Theorem~\ref{immathm:9}.
\begin{myth}\label{immathm:54.38.1}
Let $P(x,y)=a+bx+cy$ be a linear bivariate polynomial over
$\mathbb{Z}_n$. $(\mathbb{Z}_n,P)$  represents an external medial groupoid
over $\mathbb{Z}_n$ if and only if $b^2\equiv c^2(\bmod{n})$.
\end{myth}
{\bf Proof}\\
This is proved by using Lemma~\ref{immalem:50}.
\begin{myth}\label{immathm:54.38.2}
Let $P(x,y)=a+bx+cy$ be a linear bivariate polynomial over
$\mathbb{Z}_n$. $(\mathbb{Z}_n,P)$  represents an external medial quasigroup
over $\mathbb{Z}_n$ if and only if $b^2\equiv c^2(\bmod{n})$ and $(b,n)=(c,n)=1$.
\end{myth}
{\bf Proof}\\
This is proved by using Lemma~\ref{immalem:50.1} and Theorem~\ref{immathm:9}.
\begin{myexam}
$P(x,y)=4+2x+2y$ is a linear bivariate polynomial over $\mathbb{Z}_6$.
$(\mathbb{Z}_6,P)$ is an external medial groupoid over $\mathbb{Z}_6$.
\end{myexam}
\begin{myexam}
$P(x,y)=2+8x+y$ is a linear bivariate polynomial over $\mathbb{Z}_9$.
$(\mathbb{Z}_9,P)$ is an external medial quasigroup over $\mathbb{Z}_9$.
\end{myexam}
\item[(E)] Those involving left or right inverse elements:
\begin{mylem}\label{immalem:36}
Let $P(x,y)=a+bx+cy$ be a linear bivariate polynomial over
$\mathbb{Z}_n$. $P(x,y)$ represents a cross inverse property groupoid
$(\mathbb{Z}_n,P)$ over $\mathbb{Z}_n$ if and only if $a(bc-1)+x(b^2c+1-b-bc)+cy(bc-1)=0$
for all $x,y\in \mathbb{Z}_n$.
\end{mylem}
{\bf Proof}\\
$P(x,y)$ satisfies the cross inverse property $\Leftrightarrow$ $P[P(x,y),x^\rho)]=y$
$\Leftrightarrow$ $P(a+bx+cy,x^\rho)=y$ $\Leftrightarrow$ $a+b(a+bx+cy)+cx^\rho=y$ $\Leftrightarrow$ $a(bc-1)+x(b^2c+1-b-bc)+cy(bc-1)=0$ as required.
\begin{mylem}\label{immalem:36.1}
Let $P(x,y)=a+bx+cy$ be a linear bivariate polynomial over
$\mathbb{Z}_n$. $P(x,y)$ represents a cross inverse property quasigroup
$(\mathbb{Z}_n,P)$ over $\mathbb{Z}_n$ if and only if $a(bc-1)+x(b^2c+1-b-bc)+cy(bc-1)=0$
and $(b,n)=(c,n)=1$ for all $x,y,z\in \mathbb{Z}_n$.
\end{mylem}
{\bf Proof}\\
This is proved by using Lemma~\ref{immalem:36} and Theorem~\ref{immathm:9}
\begin{myth}\label{immathm:54.10}
Let $P(x,y)=a+bx+cy$ be a linear bivariate polynomial over
$\mathbb{Z}_p$ such that $a\not =0$. $P(x,y)$ represents a CIP quasigroup
$(\mathbb{Z}_p,P)$ over $\mathbb{Z}_p$ if and only if $bc\equiv 1(\bmod{p})$.
\end{myth}
{\bf Proof}\\
This is proved by using Lemma~\ref{immalem:36.1}.
\begin{myth}\label{immathm:54.11}
Let $P(x,y)=a+bx+cy$ be a linear bivariate polynomial over
$\mathbb{Z}_n$ such that $a\not =0$ and $c$ is invertible in $\mathbb{Z}_n$.
$P(x,y)$ represents a CIP groupoid
$(\mathbb{Z}_n,P)$ over $\mathbb{Z}_n$ if and only if $bc\equiv 1(\bmod{n})$.
\end{myth}
{\bf Proof}\\
This is proved by using Lemma~\ref{immalem:36}.
\begin{myth}\label{immathm:54.111}
Let $P(x,y)=a+bx+cy$ be a linear bivariate polynomial over
$\mathbb{Z}_n$ such that $a\not =0,~c$ is invertible in $\mathbb{Z}_n$ and $(b,n)=(c,n)=1$.
$P(x,y)$ represents a CIP quasigroup
$(\mathbb{Z}_n,P)$ over $\mathbb{Z}_n$ if and only if $bc\equiv 1(\bmod{n})$.
\end{myth}
{\bf Proof}\\
This is proved by using Lemma~\ref{immalem:36.1}.
\begin{myth}\label{immathm:54.12}
Let $P(x,y)=a+bx+cy$ be a linear bivariate polynomial over
$\mathbb{Z}_n$. $P(x,y)$ represents a CIP groupoid
$(\mathbb{Z}_n,P)$ over $\mathbb{Z}_n$ if $bc\equiv 1(\bmod{n})$.
\end{myth}
{\bf Proof}\\
This is proved by using Lemma~\ref{immalem:36}.
\begin{myth}\label{immathm:54.13}
Let $P(x,y)=a+bx+cy$ be a linear bivariate polynomial over
$\mathbb{Z}_n$ such that $(b,n)=(c,n)=1$. $P(x,y)$ represents a CIP quasigroup
$(\mathbb{Z}_n,P)$ over $\mathbb{Z}_n$ if $bc\equiv 1(\bmod{n})$.
\end{myth}
{\bf Proof}\\
This is proved by using Lemma~\ref{immalem:36.1}.
\begin{myexam}
$P(x,y)=2+4x+4y$ is a linear bivariate polynomial over $\mathbb{Z}_5$.
$(\mathbb{Z}_5,P)$ is a cross inverse property groupoid over $\mathbb{Z}_5$.
\end{myexam}
\begin{myexam}
$P(x,y)=3+4x+4y$ is a linear bivariate polynomial over $\mathbb{Z}_5$.
$(\mathbb{Z}_5,P)$ is a cross inverse property quasigroup over $\mathbb{Z}_5$.
\end{myexam}
\end{description}
\newpage
\begin{landscape}
\begin{table}[tbp]
\begin{center}
\begin{tabular}{|c|c|c|c|c|c|c|c|c|}
\hline
S/N & NAME & G & Q & $ \mathbb{Z}_n$ & $\mathbb{Z}_p$ & HYPO & N AND S & EXAMPLE\\
\hline
1 & Idempotent & $\checkmark$ & & $\checkmark$ & &  & $b+c=1$, $a=0$ & $5x+2y,~\mathbb{Z}_6$  \\
\hline
2 & Unipotent & $\checkmark$ & & $\checkmark$ & & & $b+c=0$ & $2+4x+2y,~\mathbb{Z}_6$\\
\cline{3-9}
 &  &  &$\checkmark$ & $\checkmark$ & & & $b+c=0$,$(b,n)=(c,n)=1$ & $2+5x+y$,~$\mathbb{Z}_6$\\
\hline
3 & Commut & $\checkmark$ & &  $\checkmark$ & & & $b=c$ & $1+4x+4y,~\mathbb{Z}_6$\\
\cline{3-9}
 &  &  & $\checkmark$ & $\checkmark$ &  & & $b=c$,$(b,n)=(c,n)=1$ & $1+5x+5y,~\mathbb{Z}_6$\\
\hline
4 & Sade Right & $\checkmark$ & & & $\checkmark$ & $a\not =0$ & $b=-1$ & $2+6x+4y,~\mathbb{Z}_7$\\
\cline{3-9}
 & & & $\checkmark$ & & $\checkmark$ & $a\not =0$ & $b=-1$ & $1+5x+4y,~\mathbb{Z}_7$\\
\hline
5 & Sade Left & $\checkmark$ & & & $\checkmark$ & $a\not =0$ & $c=-1$ & $2+4x+5y,~\mathbb{Z}_7$\\
\cline{3-9}
 & & & $\checkmark$ & & $\checkmark$ & $a\not =0$ & $c=-1$ & $2+5x+5y,~\mathbb{Z}_7$\\
\hline
6 & Right & $\checkmark$ & & & $\checkmark$ & $a\not =0$ & $b=c=1$ & $3+x+y,~\mathbb{Z}_7$\\
\cline{3-9}
 & Alternative & & $\checkmark$ & & $\checkmark$ & $a\not =0$ & $b=c=1$ & $3+x+y,~\mathbb{Z}_7$\\
\hline
7 & Left & $\checkmark$ & & & $\checkmark$ & $a\not =0$ & $b=c=1$ & $2+x+y,~\mathbb{Z}_7$\\
\cline{3-9}
 & Alternative & & $\checkmark$ & & $\checkmark$ & $a\not =0$ & $b=c=1$ & $2+x+y,~\mathbb{Z}_7$\\
\hline
8 & Medial & $\checkmark$ & & & $\checkmark$ & $a\not =0$ & $b=c$ & $2+4x+4y,~\mathbb{Z}_7$\\
\cline{3-9}
 & Alternative & $\checkmark$ & & & $\checkmark$ & $b\not =c$ & $b+c=1$ & $2+4x+2y,~\mathbb{Z}_5$\\
\cline{3-9}
 & & & $\checkmark$ & & $\checkmark$ & $a\not =0$ & $b=c$ & $2+4x+4y,~\mathbb{Z}_7$\\
\cline{3-9}
 & & & $\checkmark$ & & $\checkmark$ & $b\not =c$ & $b+c=1$ & $2+4x+2y,~\mathbb{Z}_7$\\
\hline
9 & Right & $\checkmark$ & & & $\checkmark$ & $a\not =0$ & $b=c=-1$ & $2+4x+4y,~\mathbb{Z}_5$\\
\cline{3-9}
 & Semi & $\checkmark$ & & $\checkmark$ & & $a=0$ & $bc=1,~c^2=-b$ & $5x+2y,~\mathbb{Z}_9$\\
\cline{3-9}
 & Symmetry & & $\checkmark$ & & $\checkmark$ & $a\not =0$ & $b=c=-1$ & $2+4x+4y,~\mathbb{Z}_5$\\
\cline{3-9}
 & & & $\checkmark$ & $\checkmark$ & & $a=0$ & $bc=1,~c^2=-b$ & $5x+2y,~\mathbb{Z}_9$\\
\hline
10 & Left & $\checkmark$ & & & $\checkmark$ & $a\not =0$ & $b=c=-1$ & $3+4x+4y,~\mathbb{Z}_5$\\
\cline{3-9}
 & Semi& $\checkmark$ & & $\checkmark$ & & $a=0$ & $b=1,~b^2=-c$ & $x+9y,~\mathbb{Z}_{10}$\\
\cline{3-9}
 & Symmetry & & $\checkmark$ & & $\checkmark$ & $a\not =0$ & $b=c=-1$ & $3+4x+4y,~\mathbb{Z}_5$\\
\cline{3-9}
 & & & $\checkmark$ & $\checkmark$ & & $a=0$ & $b=1,~b^2=-c$ & $x+9y,~\mathbb{Z}_{10}$\\
\hline
11 & Stein First & $\checkmark$ & & & $\checkmark$ & $a\not =0$ & $b=c$ & $3+4x+4y,~\mathbb{Z}_5$\\
\cline{3-9}
 & & & $\checkmark$ & & $\checkmark$ & $a\not =0$ & $b=c$ & $2+4x+4y,~\mathbb{Z}_5$\\
\hline
12 & Stein & $\checkmark$ & & & $\checkmark$ & $a\not =0$ & $b=c$ & $3+4x+4y,~\mathbb{Z}_5$\\
\cline{3-9}
 & Second & & $\checkmark$ & & $\checkmark$ & $a\not =0$ & $b=c$ & $2+4x+4y,~\mathbb{Z}_5$\\
\hline
13 & Schroder & $\checkmark$ & & $\checkmark$ & & & $b^2+c^2=1,~2bc=0,~a=0$ & $2x+3y,~\mathbb{Z}_6$\\
\cline{3-9}
 & Second & & $\checkmark$ & $\checkmark$ & & & $b^2+c^2=1,~2bc=0,~a=0,~(b,n)=(c,n)=1$ & ? \\
\cline{3-9}
 & & $\checkmark$ & & & $\checkmark$ & $a\not =0$ & $b+c=-1,~b^2+c^2=1,~2bc=0$ & ?\\
\cline{3-9}
 & & & $\checkmark$ & & $\checkmark$ & $a\not =0$ & $b+c=-1,~b^2+c^2=1,~2bc=0$ & ?\\
\hline
\end{tabular}
\end{center}
\end{table}
\end{landscape}
\newpage
\begin{landscape}
\begin{table}[tbp]
\begin{center}
\begin{tabular}{|c|c|c|c|c|c|c|c|c|}
\hline
14 & Stein Third & $\checkmark$ & & $\checkmark$ & & & $b^2+c^2=0,~2bc=1,~a=0$ &? \\
\cline{3-9}
 & & & $\checkmark$ & $\checkmark$ & & & $(b,n)=(c,n)=1,~b^2+c^2=0,~2bc=1,~a=0$ &? \\
\cline{3-9}
 & & $\checkmark$ & & & $\checkmark$ & $a\not =0$ & $b^2+c^2=0,~2bc=1,$ & $3+2x+4y,~\mathbb{Z}_5$\\
\cline{3-9}
 & & & $\checkmark$ & & $\checkmark$ & $a\not =0$ & $b^2+c^2=0,~2bc=1,$ & $2+2x+4y,~\mathbb{Z}_5$ \\
\hline
15 & Associative & $\checkmark$ & & & $\checkmark$ & $ a\not =0$ & $b=c=1$ & $2+x+y,~\mathbb{Z}_6$\\
\cline{3-9}
 & & & $\checkmark$ & & $\checkmark$ & $a\not =0$ & $b=c=1$ & $2+x+y,~\mathbb{Z}_6$\\
\hline
16 & Slim & $\checkmark$ & & $\checkmark$ & & $a=0,~c$ invert & $bc=0,~c=1$ &! \\
\cline{3-9}
 & & & $\checkmark$ & $\checkmark$ & & $a=0,~c$ invert & $bc=0,~c=1,~(b,n)=(c,n)=1$ &? \\
\hline
17 & Cyclic & $\checkmark$ & & $\checkmark$ & & & $b=c=1$ & $3+x+y,~\mathbb{Z}_6$\\
\cline{3-9}
 & Associativity & & $\checkmark$ & $\checkmark$ & & & $b=c=1,~(b,n)=(c,n)=1$ & $3+x+y,~\mathbb{Z}_6$\\
\hline
18 & Right & $\checkmark$ & & $\checkmark$ & & & $b=1$ & $1+x+5y,~\mathbb{Z}_6$\\
\cline{3-9}
 & Permutability & & $\checkmark$ & $\checkmark$ & & & $b=1,~(b,n)=(c,n)=1$ & $1+x+5y,~\mathbb{Z}_6$\\
\hline
19 & Left & $\checkmark$ & & $\checkmark$ & & & $c=1$ & $1+5x+y,~\mathbb{Z}_6$\\
\cline{3-9}
 & Permutability & & $\checkmark$ & $\checkmark$ & & & $c=1,~(b,n)=(c,n)=1$ & $3+5x+y,~\mathbb{Z}_6$\\
\hline
20 & Abel & $\checkmark$ & & $\checkmark$ & & &$c^2=b$ & $2+4x+2y,~\mathbb{Z}_6$\\
\cline{3-9}
 & Grassman & & $\checkmark$ & $\checkmark$ & & & $c^2=b,~(b,n)=(c,n)=1$ & $2+4x+2y,~\mathbb{Z}_9$\\
\hline
21 & Commuting & $\checkmark$ & & & $\checkmark$ & $a\not =0$ & $b=c=1$ & $1+x+y,~\mathbb{Z}_7$\\
\cline{3-9}
 & Product & & $\checkmark$ & & $\checkmark$ & $a\not =0$ & $b=c=1$ & $1+x+y,~\mathbb{Z}_7$\\
\hline
22 & Dual Comm & $\checkmark$ & & & $\checkmark$ & $a\not =0$ & $b=c=1$ & $1+x+y,~\mathbb{Z}_7$\\
\cline{3-9}
 & Product & & $\checkmark$ & & $\checkmark$ & $a\not =0$ & $b=c=1$ & $1+x+y,~\mathbb{Z}_7$\\
\hline
23 & Right & $\checkmark$ & & & $\checkmark$ & $a\not =0$ & $b=1,~c=-1$ & $2+x+6y,~\mathbb{Z}_7$\\
\cline{3-9}
 & Transitivity  & & $\checkmark$ & & $\checkmark$ & $a\not =0$ & $b=1,~c=-1$ & $2+x+6y,~\mathbb{Z}_7$\\
\hline
24 & Left & $\checkmark$ & & & $\checkmark$ & $a\not =0$ & $b=-1,~c=1$ & $2+6x+y,~\mathbb{Z}_7$\\
\cline{3-9}
 & Transitivity & & $\checkmark$ & & $\checkmark$ & $a\not =0$ & $b=-1,~c=1$ & $2+6x+y,~\mathbb{Z}_7$\\
\hline
25 & Schweitzer & $\checkmark$ & & $\checkmark$ & & $b,~c$ invert & $b=1,~c=-1$ & $2+x+5y,~\mathbb{Z}_6$\\
\cline{3-9}
 & & & $\checkmark$ & $\checkmark$ & & $b,~c$ invert & $b=1,~c=-1,~(b,n)=(c,n)=1$ & $2+x+5y,~\mathbb{Z}_6$\\
\cline{3-9}
 & & $\checkmark$ & & & $\checkmark$ & $a\not =0$ & $b=1,~c=-1$ & $3+x+6y,~\mathbb{Z}_7$\\
\cline{3-9}
 & & & $\checkmark$ & & $\checkmark$ & $a\not =0$ & $b=1,~c=-1$ & $3+x+6y,~\mathbb{Z}_7$\\
\hline
26 & Dual of & $\checkmark$ & & $\checkmark$ & & $b,~c$ invert & $b=1,~c=-1$ & $2+x+5y,~\mathbb{Z}_6$\\
\cline{3-9}
 & Schweitzer & & $\checkmark$ & $\checkmark$ & & $b,~c$ invert & $b=1,~c=-1,~(b,n)=(c,n)=1$ & $2+x+5y,~\mathbb{Z}_6$\\
\cline{3-9}
 & & $\checkmark$ & & & $\checkmark$ & $a\not =0$ & $b=1,~c=-1$ & $3+x+6y,~\mathbb{Z}_7$\\
\cline{3-9}
 & & & $\checkmark$ & & $\checkmark$ & $a\not =0$ & $b=1,~c=-1$ & $3+x+6y,~\mathbb{Z}_7$\\
\hline
\end{tabular}
\end{center}
\end{table}
\end{landscape}
\newpage
\begin{landscape}
\begin{table}[tbp]
\begin{center}
\begin{tabular}{|c|c|c|c|c|c|c|c|c|}
\hline
27 & Right Self & $\checkmark$ & & & $\checkmark$ & & $c=1-b,~a=0$ & $3x+5y,~\mathbb{Z}_7$\\
\cline{3-9}
 & Distributive & & $\checkmark$ & & $\checkmark$ & & $c=1-b,~a=0$ & $3x+5y,~\mathbb{Z}_7$\\
\hline
28 & Left Self & $\checkmark$ & & & $\checkmark$ & & $c=1-b,~a=0$ & $3x+5y,~\mathbb{Z}_7$\\
\cline{3-9}
 & Distributive & & $\checkmark$ & & $\checkmark$ & & $c=1-b,~a=0$ & $3x+5y,~\mathbb{Z}_7$\\
\hline
29 & Right & $\checkmark$ & & $\checkmark$ & & $b,~c$ invert & $b=c,~2b^2=b$ &? \\
\cline{3-9}
 & Abelian & & $\checkmark$ & $\checkmark$ & & $b,~c$ invert & $b=c,~2b^2=b$ & ?\\
\cline{3-9}
 & Distributivity & $\checkmark$ & & $\checkmark$ & & $a\not =0$ & $b=c,~2b^2=b$ & ?\\
\cline{3-9}
 & & & $\checkmark$ & $\checkmark$ & & $a\not =0$ & $b=c,~2b^2=b$ &? \\
\cline{3-9}
 & & $\checkmark$ & & & $\checkmark$ & $a\not =0$ & $b=c,~2b=1$ & $2+3x+3y,~\mathbb{Z}_5$\\
\cline{3-9}
 & & & $\checkmark$ & & $\checkmark$ & $a\not =0$ & $b=c,~2b=1$ & $2+3x+3y,~\mathbb{Z}_5$\\
\hline
30 & Left & $\checkmark$ & & $\checkmark$ & & $b,~c$ invert & $b=c,~2b^2=b$ & \\
\cline{3-9}
 & Abelian & & $\checkmark$ & $\checkmark$ & & $b,~c$ invert & $b=c,~2b^2=b$ &? \\
\cline{3-9}
 & Distributivity & $\checkmark$ & & $\checkmark$ & & $a\not =0$ & $b=c,~2b^2=b$ &? \\
\cline{3-9}
 & & & $\checkmark$ & $\checkmark$ & & $a\not =0$ & $b=c,~2b^2=b$ &? \\
\cline{3-9}
 & & $\checkmark$ & & & $\checkmark$ & $a\not =0$ & $b=c,~2b=1$ & $2+3x+3y,~\mathbb{Z}_5$\\
\cline{3-9}
 & & & $\checkmark$ & & $\checkmark$ & $a\not =0$ & $b=c,~2b=1$ & $2+3x+3y,~\mathbb{Z}_5$\\
\hline
31 & Bol & $\checkmark$ & & $\checkmark$ & & & $b=c=1$ & $2+x+y,~\mathbb{Z}_6$\\
\cline{3-9}
 & Moufang & & $\checkmark$ & $\checkmark$ & & $(b,n)=(c,n)=1$ & $b=c=1$ & $2+x+y,~\mathbb{Z}_6$\\
\hline
32 & Dual Bol & $\checkmark$ & & $\checkmark$ & & & $b=c=1$ & $2+x+y,~\mathbb{Z}_6$\\
\cline{3-9}
 & Moufang & & $\checkmark$ & $\checkmark$ & & $(b,n)=(c,n)=1$ & $b=c=1$ & $2+x+y,~\mathbb{Z}_6$\\
\hline
33 & Moufang & $\checkmark$ & & & $\checkmark$ & & $b=c=1,~a=0$ & $x+y,~\mathbb{Z}_5$\\
\cline{3-9}
 & & & $\checkmark$ & & $\checkmark$ & & $b=c=1,~a=0$ & $x+y,~\mathbb{Z}_5$\\
\hline
34 & R Bol & $\checkmark$ & & & $\checkmark$ & $a\not =0$ & $b^2=1,~b=c=1$ & $2+x+y,~\mathbb{Z}_7$\\
\cline{3-9}
 & & & $\checkmark$ & & $\checkmark$ & $a\not =0$ & $b^2=1,~b=c=1$ & $2+x+y,~\mathbb{Z}_7$\\
\cline{3-9}
 & & $\checkmark$ & & & $\checkmark$ & $-1\not =b\not =c$ & $b^2=1,c=1,a=0$ & $8x+y,~\mathbb{Z}_{63}$\\
\cline{3-9}
 & & & $\checkmark$ & & $\checkmark$ & $-1\not =b\not =c$ & $b^2=1,c=1,a=0$ & $8x+y,~\mathbb{Z}_{63}$\\
\hline
35 & L Bol & $\checkmark$ & & & $\checkmark$ & $a\not =0$ & $c^2=1,~b=c=1$ & $2+x+y,~\mathbb{Z}_7$\\
\cline{3-9}
 & & & $\checkmark$ & & $\checkmark$ & $a\not =0$ & $c^2=1,~b=c=1$ & $2+x+y,~\mathbb{Z}_7$\\
\cline{3-9}
 & & $\checkmark$ & & & $\checkmark$ & $-1\not =b\not =c$ & $c^2=1,b=1,a=0$ & $x+8y,~\mathbb{Z}_{63}$\\
\cline{3-9}
 & & & $\checkmark$ & & $\checkmark$ & $-1\not =b\not =c$ & $c^2=1,b=1,a=0$ & $x+8y,~\mathbb{Z}_{63}$\\
\hline
\end{tabular}
\end{center}
\end{table}
\end{landscape}
\newpage
\begin{landscape}
\begin{table}[tbp]
\begin{center}
\begin{tabular}{|c|c|c|c|c|c|c|c|c|}
\hline
36 & $\textrm{RC}_\textrm{4}$ & $\checkmark$ & & & $\checkmark$ & $a=0$ & $c=b^2=1$ & $8x+y,~\mathbb{Z}_{63}$\\
\cline{3-9}
 & & & $\checkmark$ & & $\checkmark$ & $a=0$ & $c=b^2=1$ & $8x+y,~\mathbb{Z}_{63}$\\
\cline{3-9}
 & & $\checkmark$ & & $\checkmark$ & & $a=0,~b,c$ invert & $c=b^2=1$ & $8x+y,~\mathbb{Z}_{63}$\\
\cline{3-9}
 & & & $\checkmark$ & $\checkmark$ & & $a=0,~b,c$ invert & $c=b^2=1,~(b,n)=(c,n)=1$ & $8x+y,~\mathbb{Z}_{63}$\\
\cline{3-9}
 & & $\checkmark$ & & $\checkmark$ & & & $b=-1,~c=1$ & $2+5x+y,~\mathbb{Z}_6$\\
\cline{3-9}
 & & & $\checkmark$ & $\checkmark$ & & & $b=-1,~c=1,~(b,n)=(c,n)=1$ & $2+5x+y,~\mathbb{Z}_6$\\
\hline
37 & $\textrm{LC}_\textrm{4}$ & $\checkmark$ & & & $\checkmark$ & $a=0$ & $b=c^2=1$ & $x+8y,~\mathbb{Z}_{63}$\\
\cline{3-9}
 & & & $\checkmark$ & & $\checkmark$ & $a=0$ & $b=c^2=1$ & $x+8y,~\mathbb{Z}_{63}$\\
\cline{3-9}
 & & $\checkmark$ & & $\checkmark$ & & $a=0,~b,c$ invert & $b=c^2=1$ & $x+3y,~\mathbb{Z}_8$\\
\cline{3-9}
 & & & $\checkmark$ & $\checkmark$ & & $a=0,~b,c$ invert & $b=c^2=1,~(b,n)=(c,n)=1$ & $x+4y,~\mathbb{Z}_{15}$\\
\cline{3-9}
 & & $\checkmark$ & & $\checkmark$ & & & $b=-1,~c=1$ & $2+5x+y,~\mathbb{Z}_6$\\
\cline{3-9}
 & & & $\checkmark$ & $\checkmark$ & & & $b=-1,~c=1,~(b,n)=(c,n)=1$ & $2+5x+y,~\mathbb{Z}_6$\\
\hline
38 & $\textrm{RC}_\textrm{1}$ & $\checkmark$ & & & $\checkmark$ & $a=0$ & $c=b^2=1$ & $8x+y,~\mathbb{Z}_{63}$\\
\cline{3-9}
 & & & $\checkmark$ & & $\checkmark$ & $a=0$ & $c=b^2=1$ & $8x+y,~\mathbb{Z}_{63}$\\
\cline{3-9}
 & & $\checkmark$ & & $\checkmark$ & & $a=0,~b,c$ invert & $c=b^2=1$ & $8x+y,~\mathbb{Z}_{63}$\\
\cline{3-9}
 & & & $\checkmark$ & $\checkmark$ & & $a=0,~b,c$ invert & $c=b^2=1,~(b,n)=(c,n)=1$ & $8x+y,~\mathbb{Z}_{63}$\\
\cline{3-9}
 & & $\checkmark$ & & $\checkmark$ & & & $b=-1,~c=1$ & $2+5x+y,~\mathbb{Z}_6$\\
\cline{3-9}
 & & & $\checkmark$ & $\checkmark$ & & & $b=-1,~c=1,~(b,n)=(c,n)=1$ & $2+5x+y,~\mathbb{Z}_6$\\
\hline
39 & $\textrm{LC}_\textrm{1}$ & $\checkmark$ & & & $\checkmark$ & $a=0,~c\not =1$ & $c=-1$ & $3x+6y,~\mathbb{Z}_7$\\
\cline{3-9}
 & & & $\checkmark$ & & $\checkmark$ & $a=0,~c\not =1$ & $c=-1$ & $3x+6y,~\mathbb{Z}_7$\\
\cline{3-9}
 & & $\checkmark$ & & $\checkmark$ & & $a=0,~c\not =1,~c$ invert & $c=-1$ & $5x+5y,~\mathbb{Z}_6$\\
\cline{3-9}
 & & & $\checkmark$ & $\checkmark$ & & $a=0,~c\not =1,~c$ invert & $c=-1,~(b,n)=(c,n)=1$ & $5x+5y,~\mathbb{Z}_6$\\
\hline
40 & $\textrm{LC}_\textrm{3}$ & $\checkmark$ & & $\checkmark$ & & & $c=1,~b=-2$ & $3+4x+y,~\mathbb{Z}_6$\\
\cline{3-9}
 & & & $\checkmark$ & $\checkmark$ & & & $c=1,~b=-2,~(b,n)=(c,n)=1$ & $2+5x+y,~\mathbb{Z}_7$\\
\hline
41 & $\textrm{RC}_\textrm{3}$ & $\checkmark$ & & $\checkmark$ & & & $c=1,~b=-2$ & $3+4x+y,~\mathbb{Z}_6$\\
\cline{3-9}
 & & & $\checkmark$ & $\checkmark$ & & & $c=1,~b=-2,~(b,n)=(c,n)=1$ & $2+5x+y,~\mathbb{Z}_7$\\
\hline
42 & C-Law & $\checkmark$ & & & $\checkmark$ & $a=0$ & $b=c=-1$ & $4x+4y,~\mathbb{Z}_5$\\
\cline{3-9}
 & & & $\checkmark$ & & $\checkmark$ & $a=0$ & $b=c=-1$ & $4x+4y,~\mathbb{Z}_5$\\
\cline{3-9}
 & & $\checkmark$ & & $\checkmark$ & & $a\not =0,~b\not =1,~b,c$ inv & $b=c=-1$ & $3+5x+5y,~\mathbb{Z}_6$\\
\cline{3-9}
 & & & $\checkmark$ & $\checkmark$ & & $a\not =0,~b\not =1,~b,c$ inv & $b=c=-1,~(b,n)=(c,n)=1$ & $3+5x+5y,~\mathbb{Z}_6$\\
 \hline
43 & LIP & $\checkmark$ & & & $\checkmark$ & $a\not =0$ & $c^2=b^2=bc=1$ &?\\
\cline{3-9}
 & & & $\checkmark$ & & $\checkmark$ & $a\not =0$ & $c^2=b^2=bc=1$ &? \\
 \hline
44 & RIP & $\checkmark$ & & & $\checkmark$ & $a\not =0$ & $c^2=b^2=bc=1$ &? \\
\cline{3-9}
 & & & $\checkmark$ & & $\checkmark$ & $a\not =0$ & $c^2=b^2=bc=1$ &? \\
\hline
\end{tabular}
\end{center}
\end{table}
\end{landscape}
\newpage
\begin{landscape}
\begin{table}[tbp]
\begin{center}
\begin{tabular}{|c|c|c|c|c|c|c|c|c|}
\hline
45 & 1st Right & $\checkmark$ & & & $\checkmark$ & $a\not =0$ & $bc=1$ & $2+3x+4y,~\mathbb{Z}_{11}$\\
\cline{3-9}
 & CIP & & $\checkmark$ & & $\checkmark$ & $a\not =0$ & $bc=1$ & $2+3x+4y,~\mathbb{Z}_{11}$\\
\cline{3-9}
 & & $\checkmark$ & & $\checkmark$ & & $a\not =0,~c$ inv & $bc=1$ & $3+3x+3y,~\mathbb{Z}_8$\\
\cline{3-9}
 & & & $\checkmark$ & $\checkmark$ & & $a\not =0,~c$ inv & $bc=1,~(b,n)=(c,n)=1$ & $3+3x+3y,~\mathbb{Z}_8$\\
\hline
46 & 2nd Right & $\checkmark$ & & $\checkmark$ & & & $bc=1$ & $3+3x+3y,~\mathbb{Z}_8$\\
\cline{3-9}
 & CIP & & $\checkmark$ & $\checkmark$ & & & $bc=1,~(b,n)=(c,n)=1$ & $3+3x+3y,~\mathbb{Z}_8$\\
\hline
47 & 1st Left & $\checkmark$ & & & $\checkmark$ & $a\not =0$ & $bc=1$ & $2+3x+4y,~\mathbb{Z}_{11}$\\
\cline{3-9}
 & CIP & & $\checkmark$ & & $\checkmark$ & $a\not =0$ & $bc=1$ & $2+3x+4y,~\mathbb{Z}_{11}$\\
\cline{3-9}
 & & $\checkmark$ & & $\checkmark$ & & $a\not =0,~b$ inv & $bc=1$ & $3+3x+3y,~\mathbb{Z}_8$\\
\cline{3-9}
 & & & $\checkmark$ & $\checkmark$ & & $a\not =0,~b$ inv & $bc=1,~(b,n)=(c,n)=1$ & $3+3x+3y,~\mathbb{Z}_8$\\
\hline
48 & 2nd Left & $\checkmark$ & & $\checkmark$ & & & $bc=1$ & $3+3x+3y,~\mathbb{Z}_8$\\
\cline{3-9}
 & CIP & & $\checkmark$ & $\checkmark$ & & & $bc=1,~(b,n)=(c,n)=1$ & $3+3x+3y,~\mathbb{Z}_8$\\
\hline
49 & R AAIP & $\checkmark$ & & & $\checkmark$ & $bc+b\not =1$ & $b=c$ & $2+4x+4y,~\mathbb{Z}_{11}$\\
\cline{3-9}
 & & & $\checkmark$ & & $\checkmark$ & $bc+b\not =1$ & $b=c$ & $2+4x+4y,~\mathbb{Z}_{11}$\\
\cline{3-9}
 & & $\checkmark$ & & & $\checkmark$ & $c\not =b$ & $b+bc=1$ & $2+3x+y,~\mathbb{Z}_5$\\
\cline{3-9}
 & & & $\checkmark$ & & $\checkmark$ & $c\not =b$ & $b+bc=1$ & $2+3x+y,~\mathbb{Z}_5$\\
\hline
50 & L AAIP & $\checkmark$ & & & $\checkmark$ & $bc+b\not =1$ & $b=c$ & $2+4x+4y,~\mathbb{Z}_{11}$\\
\cline{3-9}
 & & & $\checkmark$ & & $\checkmark$ & $bc+b\not =1$ & $b=c$ & $2+4x+4y,~\mathbb{Z}_{11}$\\
\cline{3-9}
 & & $\checkmark$ & & & $\checkmark$ & $c\not =b$ & $b+bc=1$ & $2+3x+y,~\mathbb{Z}_5$\\
\cline{3-9}
 & & & $\checkmark$ & & $\checkmark$ & $c\not =b$ & $b+bc=1$ & $2+3x+y,~\mathbb{Z}_5$\\
\hline
51 & R AIP & $\checkmark$ & & $\checkmark$ & & & & $a+bx+cy,~\mathbb{Z}_n$\\
\cline{3-9}
 & & & $\checkmark$ & $\checkmark$ & & & & $a+bx+cy,~\mathbb{Z}_n$\\
\cline{3-9}
 & & $\checkmark$ & & & $\checkmark$ & & & $a+bx+cy,~\mathbb{Z}_n$\\
\cline{3-9}
 & & & $\checkmark$ & & $\checkmark$ & & & $a+bx+cy,~\mathbb{Z}_n$\\
\hline
52 & L AIP & $\checkmark$ & & $\checkmark$ & & & & $a+bx+cy,~\mathbb{Z}_n$\\
\cline{3-9}
 & & & $\checkmark$ & $\checkmark$ & & & & $a+bx+cy,~\mathbb{Z}_n$\\
\cline{3-9}
 & & $\checkmark$ & & & $\checkmark$ & & & $a+bx+cy,~\mathbb{Z}_n$\\
\cline{3-9}
 & & & $\checkmark$ & & $\checkmark$ & & & $a+bx+cy,~\mathbb{Z}_n$\\
\hline
53 & R SAIP & $\checkmark$ & & $\checkmark$ & & & & $a+bx+cy,~\mathbb{Z}_n$\\
\cline{3-9}
 & & & $\checkmark$ & $\checkmark$ & & & & $a+bx+cy,~\mathbb{Z}_n$\\
\cline{3-9}
 & & $\checkmark$ & & & $\checkmark$ & & & $a+bx+cy,~\mathbb{Z}_n$\\
\cline{3-9}
 & & & $\checkmark$ & & $\checkmark$ & & & $a+bx+cy,~\mathbb{Z}_n$\\
\hline
\end{tabular}
\end{center}
\end{table}
\end{landscape}
\newpage
\begin{landscape}
\begin{table}[!hbp]
\begin{center}
\begin{tabular}{|c|c|c|c|c|c|c|c|c|}
\hline
54 & L SAIP & $\checkmark$ & & $\checkmark$ & & & & $a+bx+cy,~\mathbb{Z}_n$\\
\cline{3-9}
 & & & $\checkmark$ & $\checkmark$ & & & & $a+bx+cy,~\mathbb{Z}_n$\\
\cline{3-9}
 & & $\checkmark$ & & & $\checkmark$ & & & $a+bx+cy,~\mathbb{Z}_n$\\
\cline{3-9}
 & & & $\checkmark$ & & $\checkmark$ & & & $a+bx+cy,~\mathbb{Z}_n$\\
\hline
55 & R WIP & $\checkmark$ & & & $\checkmark$ & $a=0,~c^2\not =0$ & $bc=1$ & $3x+5y,~\mathbb{Z}_7$\\
\cline{3-9}
 & & & $\checkmark$ & & $\checkmark$ & $a=0,~c^2\not =0$ & $bc=1$ & $3x+5y,~\mathbb{Z}_7$\\
\cline{3-9}
 & & $\checkmark$ & & $\checkmark$ & & $a=0,~c$ inv & $bc=1$ & $3x+4y,~\mathbb{Z}_6$\\
\cline{3-9}
 & & & $\checkmark$ & $\checkmark$ & & $a=0,~c$ inv & $bc=1,~(b,n)=(c,n)=1$ &? \\
\cline{3-9}
 & & $\checkmark$ & & $\checkmark$ & & $a=0,~bc+b\not =1$ & $bc=1$ & ?\\
\cline{3-9}
 & & & $\checkmark$ & $\checkmark$ & & $a=0,~bc+b\not =1$ & $bc=1,~(b,n)=(c,n)=1$ & ?\\
\hline
56 & L WIP & $\checkmark$ & & & $\checkmark$ & $a=0,~b^2\not =0$ & $bc=1$ & $3x+5y,~\mathbb{Z}_7$\\
\cline{3-9}
 & & & $\checkmark$ & & $\checkmark$ & $a=0,~b^2\not =0$ & $bc=1$ & $3x+5y,~\mathbb{Z}_7$\\
\cline{3-9}
 & & $\checkmark$ & & $\checkmark$ & & $a=0,~b$ inv & $bc=1$ & $3x+4y,~\mathbb{Z}_6$\\
\cline{3-9}
 & & & $\checkmark$ & $\checkmark$ & & $a=0,~b$ inv & $bc=1,~(b,n)=(c,n)=1$ &? \\
\cline{3-9}
 & & $\checkmark$ & & $\checkmark$ & & $a=0,~bc+c\not =1$ & $bc=1$ & ?\\
\cline{3-9}
 & & & $\checkmark$ & $\checkmark$ & & $a=0,~bc+c\not =1$ & $bc=1,~(b,n)=(c,n)=1$ &? \\
\hline
57 & $\textrm{E}_\textrm{l}$ & $\checkmark$ & $\checkmark$ & & & & & $a+bx+cy,~\mathbb{Z}_n$\\
\cline{3-9}
 & & & $\checkmark$ & $\checkmark$ & & & & $a+bx+cy,~\mathbb{Z}_n$\\
\cline{3-9}
 & & $\checkmark$ & & & $\checkmark$ & & & $a+bx+cy,~\mathbb{Z}_n$\\
\cline{3-9}
 & & & $\checkmark$ & & $\checkmark$ & & & $a+bx+cy,~\mathbb{Z}_n$\\
\hline
58 & $\textrm{E}_\textrm{r}$ & $\checkmark$ & $\checkmark$ & & & & & $a+bx+cy,~\mathbb{Z}_n$\\
\cline{3-9}
 & & & $\checkmark$ & $\checkmark$ & & & & $a+bx+cy,~\mathbb{Z}_n$\\
\cline{3-9}
 & & $\checkmark$ & & & $\checkmark$ & & & $a+bx+cy,~\mathbb{Z}_n$\\
\cline{3-9}
 & & & $\checkmark$ & & $\checkmark$ & & & $a+bx+cy,~\mathbb{Z}_n$\\
\hline
59 & Right F & $\checkmark$ & $\checkmark$ & & & & & $a+bx+cy,~\mathbb{Z}_n$\\
\cline{3-9}
 & & & $\checkmark$ & $\checkmark$ & & & & $a+bx+cy,~\mathbb{Z}_n$\\
\cline{3-9}
 & & $\checkmark$ & & & $\checkmark$ & & & $a+bx+cy,~\mathbb{Z}_n$\\
\cline{3-9}
 & & & $\checkmark$ & & $\checkmark$ & & & $a+bx+cy,~\mathbb{Z}_n$\\
\hline
60 & Left F & $\checkmark$ & $\checkmark$ & & & & & $a+bx+cy,~\mathbb{Z}_n$\\
\cline{3-9}
 & & & $\checkmark$ & $\checkmark$ & & & & $a+bx+cy,~\mathbb{Z}_n$\\
\cline{3-9}
 & & $\checkmark$ & & & $\checkmark$ & & & $a+bx+cy,~\mathbb{Z}_n$\\
\cline{3-9}
 & & & $\checkmark$ & & $\checkmark$ & & & $a+bx+cy,~\mathbb{Z}_n$\\
\hline
\end{tabular}
\end{center}
\end{table}
\end{landscape}
\newpage
\begin{table}[!hbp]
\begin{center}
\begin{tabular}{|c|c|c|c|c|c|c|c|c|}
\hline
61 & Medial & $\checkmark$ & $\checkmark$ & & & & & $a+bx+cy,~\mathbb{Z}_n$\\
\cline{3-9}
 & & & $\checkmark$ & $\checkmark$ & & & & $a+bx+cy,~\mathbb{Z}_n$\\
\cline{3-9}
 & & $\checkmark$ & & & $\checkmark$ & & & $a+bx+cy,~\mathbb{Z}_n$\\
\cline{3-9}
 & & & $\checkmark$ & & $\checkmark$ & & & $a+bx+cy,~\mathbb{Z}_n$\\
\hline
62 & Specialized & $\checkmark$ & $\checkmark$ & & & & & $a+bx+cy,~\mathbb{Z}_n$\\
\cline{3-9}
 & Medial & & $\checkmark$ & $\checkmark$ & & & & $a+bx+cy,~\mathbb{Z}_n$\\
\cline{3-9}
 & & $\checkmark$ & & & $\checkmark$ & & & $a+bx+cy,~\mathbb{Z}_n$\\
\cline{3-9}
 & & & $\checkmark$ & & $\checkmark$ & & & $a+bx+cy,~\mathbb{Z}_n$\\
\hline
63 & First & $\checkmark$ & & & $\checkmark$ & & $b=c$ & $2+4x+4y,~\mathbb{Z}_7$\\
\cline{3-9}
 & Rectangle & & $\checkmark$ & & $\checkmark$ & & $b=c$ & $2+4x+4y,~\mathbb{Z}_7$\\
\cline{3-9}
 & & $\checkmark$ & & $\checkmark$ & & $c$ inv & $b=c$ & $2+4x+4y,~\mathbb{Z}_6$\\
\cline{3-9}
 & & & $\checkmark$ & $\checkmark$ & & $c$ inv & $b=c,~(b,n)=(c,n)=1$ & $2+4x+4y,~\mathbb{Z}_6$\\
\hline
64 & Second & $\checkmark$ & & & $\checkmark$ & & $b=-c$ & $2+4x+4y,~\mathbb{Z}_7$\\
\cline{3-9}
 & Rectangle & & $\checkmark$ & & $\checkmark$ & & $b=-c$ & $2+4x+4y,~\mathbb{Z}_7$\\
\cline{3-9}
 & & $\checkmark$ & & $\checkmark$ & & $b$ inv & $b=-c$ & $2+4x+4y,~\mathbb{Z}_6$\\
\cline{3-9}
 & & & $\checkmark$ & $\checkmark$ & & $b$ inv & $b=-c,~(b,n)=(c,n)=1$ & $2+4x+4y,~\mathbb{Z}_6$\\
\hline
65 & $\textrm{C}_i,~i=1-6$ & & $\checkmark$ & & $\checkmark$ & & $b=c$ & $3+5x+5y,~\mathbb{Z}_7$\\
\hline
66 & $\textrm{CM}_i,~i=1-14$ & & $\checkmark$ & & $\checkmark$ & $b\not =-c$ & $b=c$ & $3+5x+5y,~\mathbb{Z}_7$\\
\hline
\end{tabular}
\end{center}\caption{A Table for the Characterization of Varieties of Groupoids and Quasigroups Generated by $P(x,y)$ over $\mathbb{Z}_n$}\label{osbornloop9}
\end{table}
\begin{myrem}
A summary of the results on the characterization of groupoids and quasigroups generated by $P(x,y)$ is exhibited in Table~\ref{osbornloop9}. In the table "$G$" stands for "groupoid", "$Q$" stands for "quasigroup", "HYPO" stands for "hypothesis", "N AND S" stands for "necessary and sufficient condition(s)". Cells with question marks mean examples could not be gotten.
\end{myrem}

\end{document}